\documentclass[12pt,a4paper]{article}
\usepackage[width=18cm,height=23cm]{geometry}
\usepackage{amsmath,amssymb,amsfonts,amsthm,undertilde,enumerate}
\providecommand{\keywords}[1]{\textbf{\textit{Keywords: }}#1}

\usepackage{graphicx}
\usepackage{tikz}
\usetikzlibrary{shapes,backgrounds}
\usepackage{tikz-cd}
\usetikzlibrary{%
  matrix,%
  calc,%
  arrows%
}

\DeclareMathOperator{\im}{im}
\DeclareMathOperator{\res}{\text{res}}  
\usepackage[small,nohug,heads=vee]{diagrams}
\diagramstyle[labelstyle=\scriptstyle]
\usepackage{mathrsfs,amsbsy,bm,mathtools}
\usepackage{indent first}

\newcommand{\comm}{\circlearrowright}
\newcommand{\bc}{\mathbb{C}}
\newcommand{\cc}{\mathscr{C}}
\newcommand{\bp}{\mathbb{P}}
\newcommand{\bz}{\mathbb{Z}}
\newcommand{\br}{\mathbb{R}}
\newcommand{\co}{\mathscr{O}}
\newcommand{\cf}{\mathscr{F}}
\newcommand{\fglr}{\mathcal{GL}_r(\co)}
\newcommand{\cu}{\mathscr{U}}
\newcommand{\al}[1]{\begin{align*}{#1}\end{align*}}
\newcommand{\cond}[1]{\quad{\scriptstyle{(#1)}}}
\newcommand{\rw}{\longrightarrow}

\newcommand{\xrw}{\xrightarrow}

\newtheorem{Thm}{Theorem}[section]%
\newtheorem{Lem}[Thm]{Lemma}%
\newtheorem{Cor}[Thm]{Corollary}%
\newtheorem{Prop}[Thm]{Proposition}%
\newtheorem{Rmk}[Thm]{Remark}%
\newtheorem{Quest}[Thm]{Question}%

\title{\bf A Counterexample to Hartogs' Type Extension\\ of Holomorphic Line Bundles}
\author{Zhangchi Chen \footnote{zhangchi.chen@u-psud.fr} \footnote{D\'epartement de Math\'ematiques D'orsay, Facult\'e des Sciences, Universit\'e Paris-Sud, B\^atiment 425, 91405 Orsay Cedex, France}}
\date{\today}

\begin{document}
\maketitle

\abstract{Consider a domain $\varOmega$ in $\bc^n$ with $n\geqslant 2$ and a compact subset $K\subset\varOmega$ such that $\varOmega\backslash K$ is connected. We address the problem whether a holomorphic line bundle defined on $\varOmega\backslash K$ extends to $\varOmega$. In 2013, Forn\ae ss, Sibony and Wold gave a positive answer in dimension $n\geqslant 3$, when $\varOmega$ is pseudoconvex and $K$ is a sublevel set of a strongly plurisubharmonic exhaustion function. However, for $K$ of general shape, we construct counterexamples in any dimension $n\geqslant 2$. The key is a certain gluing lemma by means of which we extend any two holomorphic line bundles which are isomorphic on the intersection of their base spaces.}





\section{Introduction}
\label{sect-intro}

The Hartogs' extension theorem is one of the most distinctive results
in several complex variables. Let $\varOmega\subset \bc^n$
{\small($n\geqslant 2$)} be a domain. Let $K\subset\subset \varOmega$ be
a compact subset such that $\varOmega\backslash K$ is connected. Denote
by $\co$ the sheaf of holomorphic functions on $\bc^n$.

\begin{Thm}\label{Hart-1} {\em (Hartogs' extension theorem for holomorphic functions)} The restriction map
\begin{align*}
H^0(\varOmega,\co)\rw H^0(\varOmega\backslash K,\co)
\end{align*}
is bijective.
\end{Thm}

A proof using no $\overline{\partial}$ techniques
can be found in Merker-Porten's paper \cite{JM-P}.

Next, let $\co^*$ be the sheaf of invertible holomorphic functions 
on $\bc^n$. Arguing that a nowhere vanishing function $f \in 
\varOmega \backslash K$ extends holomorphically to $\varOmega$ as well
as its inverse $g := \frac{1}{f}$, and that $f \, g \equiv 1$
transfers from $\varOmega \backslash K$ to $\varOmega$
by the uniqueness principle, one deduces the

\begin{Cor}\label{Hart-2} {\em (Extension of invertible holomorphic functions)} The restriction map
\begin{align}\label{res-inv}
H^0(\varOmega,\co^*)\rw H^0(\varOmega\backslash K,\co^*)
\end{align}
is bijective.\qed
\end{Cor}

Beyond functions, it is natural to ask whether 
for holomorphic line bundles, 
Hartogs' type extension holds from $\varOmega\backslash K$ to $\varOmega$. 
If yes, is the extension unique modulo
isomorphism? 

Recall that there is a bijection between the set of isomorphic classes
of holomorphic line bundles over $\varOmega$, and the Picard group
$H^1(\varOmega,\co^*)$, constructed in the following way. Any holomorphic
line bundle $\pi \colon L \longrightarrow \varOmega$ admits an open cover
$\{U_i\}$ of $\varOmega$ together with local trivialization maps
$\varphi_i \colon \pi^{-1}(U_i) \overset{\sim}{\longrightarrow} U_i
\times \mathbb{C}$ and transition maps $f_{ij}\in H^0(U_i\cap
U_j,\co^*)$. The data $\{f_{ij}\}$ is a \v{C}ech 1-cocycle
representing some element in $H^1(\{U_i\},\co^*)\hookrightarrow
H^1(\varOmega,\co^*).$ Reciprocally, any element in $H^1(\varOmega,\co^*)$
can be expressed by some \v{C}ech 1-cocycle $\{f_{ij}\}$ with respect
to some open cover $\{U_i\}$ of $\varOmega$ valued in $\co^*$. The data
($\{U_i\}$, $\{f_{ij}\}$) gives a holomorphic line bundle.

Using these notations, we may restate our question more
precisely.

\begin{Quest}\label{Quest-ext-L}
Given a holomorphic line bundle $L$ over $\varOmega\backslash K$, does there exist a holomorphic line bundle $\widetilde{L}$ over $\varOmega$ such that $\widetilde{L}|_{\varOmega\backslash K}\cong L$? Equivalently, is the restriction 
map
\begin{align}\label{H-1-O-star-Omega-K}
H^1\big(
\varOmega,\co^*
\big)
\longrightarrow
H^1\big(
\varOmega\backslash K,\co^*
\big)
\end{align}
surjective? If yes, is it bijective?
\end{Quest}

A positive answer, under certain circumstances, was given by Fornaess-Sibony-Wold in \cite{FSW}.

\begin{Thm}\label{FSW} {\em (Extension across strictly pseudoconcave level sets)} Let $\varOmega\subset\bc^n$ {\small($n\geqslant 3$)} be a pseudoconvex domain with a $\cc^{\infty}$ strictly plurisubharmonic (psh) exhaustion function $\rho$, i.e. for each $a\in\br$, the sublevel set $K_a:=\rho^{-1}(-\infty,a]$ is compact in $\varOmega$. Then every holomorphic line bundle over $\varOmega\backslash K_a$ extends to $\varOmega$. The extension is unique modulo isomorphism.
\end{Thm}

Actually they proved a stronger version of this theorem, namely existence (resp. uniqueness) of an extension when the Levi form of $\rho$ has at least 3 (resp. 2) positive eigenvalues.

The proof of Theorem~\ref{FSW} uses (1) the exponential sequence and Cartan's theorem B (2) the extension of holomorphic functions across a totally real plane and (3) Andreotti-Grauert theory. We will present the first two ingredients in Section~\ref{sec-back} because we are going to use them later.

Now, let us come back to Question~\ref{Quest-ext-L}. For $n = 2$, Ivashkovich already presented in~\cite{Iva} a local counterexample (cex), but with $K \subset \varOmega$ not compact. In Section~\ref{sect-counter-dim-2}, we will briefly restate his construction, and by taking exponential, we will produce a domain $\varOmega \subset \mathbb{C}^2$ and a compact $K \subset\subset \varOmega$ through which some holomorphic line bundles do not extend.

\begin{Prop}
There exists a bounded pseudoconvex domain $\varOmega \subset 
\mathbb{C}^2$ equipped with a $\mathscr{C}^\infty$ strictly 
psh function
\[
\rho\colon\ \ \
\varOmega
\,\longrightarrow\,
[0,\infty)
\]
such that $K := \rho^{-1}(0) \cong S^1 \times S^1$ is a compact
totally real $2$-torus, and there exists a (nontrivial) holomorphic
line bundle $L$ on $\varOmega \backslash K$ having the property
that there exists {\em no} holomorphic line bundle
$\widetilde{L}$ on $\varOmega$ with $\widetilde{L}
\big\vert_{\varOmega \backslash K}\cong L$.
\end{Prop}

However, a similar construction in dimension $n \geqslant 3$, again
with a compact $K = \rho^{-1}(0) \cong (S^1)^n$ of the
same kind, would fall under
the {\em positive} (known) extension Theorem~\ref{FSW}.

Hence to really produce a cex to Hartogs' type extension 
for holomorphic line bundles in {\em all} dimensions $n \geqslant
2$, the compact $K\subset\subset \varOmega$ should not be
of the shape $\{\rho\leqslant a\}$, i.e. a sublevel set of a 
strictly psh exhaustion function.

In Section~\ref{sect-counter-dim-3}, we will perform
an alternative construction. In $\mathbb{C}^n$
{\small ($n \geqslant 2$)}, 
for $0 < \epsilon < n$, we introduce the domain:
\[
G_\epsilon
\,:=\,
\big\{
z\in\mathbb{C}^n
\colon\,
\sum_{j=1}^n\,
\big(\log\,\vert z_j\vert\big)^2
<
\epsilon
\big\},
\]
which contains the $n$-dimensional standard totally real torus:
\[
\mathbb{T}^n
\,=\,
\big\{
\vert z_1\vert
=\cdots=
\vert z_n\vert=1
\big\}
\,\cong\,
(S^1)^n.
\]
For $0 < \epsilon \ll n$ small, $G_\epsilon$ will appear to be
a thin Grauert tube around $\mathbb{T}^n$. 
We will check that the domain $G_\epsilon$ is relatively
compact in the ball:
\[
\varOmega
\,:=\,
B
\big(
2\sqrt{n}\,e^{\sqrt{\epsilon}}
\big)
\]
centered at the origin and of radius $2\sqrt{n}\, e^{ \sqrt{\epsilon}}
$. Also, we will take a small open ball $U_p \subset
\mathbb{C}^n$ centered at the point:
\[
p
\,:=\,
\big(
e^{\sqrt{\epsilon/n}},
\dots,
e^{\sqrt{\epsilon/n}}
\big)
\,\in\,
\partial
G_\epsilon,
\]
as connecting the interior and the exterior of $\partial G_\varepsilon$ through a small hole at $p$. Our main result is the

\begin{Thm}
\label{main-counterexample}
With the compact:
\[
K
\,:=\,
\partial G_\epsilon
\big\backslash
U_p,
\]
the open set $\varOmega \backslash K$ is connected, and there 
exists a (nontrivial) holomorphic
line bundle $L_{\text{cex}}$ on $\varOmega \backslash K$ having the property
that there exists {\em no} holomorphic line bundle
$\widetilde{L}$ on $\varOmega$ with $\widetilde{L}
\big\vert_{\varOmega \backslash K}\cong L_{\text{cex}}$. Here `cex' stands for `counterexample'
\end{Thm}

The way we construct this non-extendable $L_{\text{cex}}$ is by using the following gluing lemma.
\begin{Lem}\label{gluing-lemma}
Let $U,V\subset \bc^n$ be two open subsets, $L_U,L_V$ be two holomorphic line bundles defined over $U$ and $V$ respectively. If $L_U|_{U\cap V}\cong L_V|_{V\cap U}$ are isomorphic as holomorphic line bundles, then there exists a holomorphic line bundle $L$ defined over $U\cup V$ such that $L|_U\cong L_U,L|_V\cong L_V$.
\end{Lem}

A more general version of this gluing lemma, for holomorphic vector bundles, is stated and proved in subsection \ref{subsection-gluing}.

Note that in this lemma, we assume no geometrical condition on $U,V$ and no triviality of $L_U,L_V$. The only condition is that $L_U|_{U\cap V}\cong L_V|_{V\cap U}$. In particular, when $H^1(U\cap V,\co^*)=0$, e.g. when $U\cap V$ is convex, this condition is always satisfied.

The picard group $H^1(G_\epsilon,\co^*)\cong \bz^{\binom{n}{2}}$ is nontrivial, which will be proved in Proposition \ref{prop:G}. So we can take a nontrivial holomorphic line bundle $L_{\text{nt}}$ over $G_\epsilon$. As a consequence of Proposition \ref{real-hessian} we show there exists a small ball $U_p$ centered at $p\in\partial G_\epsilon$ such that $U_p\cap G_\epsilon$ is convex. So in the gluing lemma, if we regard $U$ as $G_\epsilon$ and $V$ as $(\varOmega\backslash\bar{U})\cup U_p$, then $U\cap V=U_p\cap G_\epsilon$ is convex. Thus we can glue $L_{\text{nt}}$ with a trivial line bundle $L_{\text{triv}}$ over $V$ to obtain a line bundle $L_{\text{cex}}$ over $U\cup V=\varOmega\backslash(\partial G_\epsilon\backslash U_p)=\varOmega\backslash K$, which is connected by Proposition \ref{prop-omega}. Such $L_{\text{cex}}$ is nontrivial since $L|_U\cong L_{\text{nt}}$ is. It cannot be extended to $\varOmega$ since $H^1(\varOmega,\co^*)=0$.
\begin{center}
\begin{picture}(0,0)%
\includegraphics{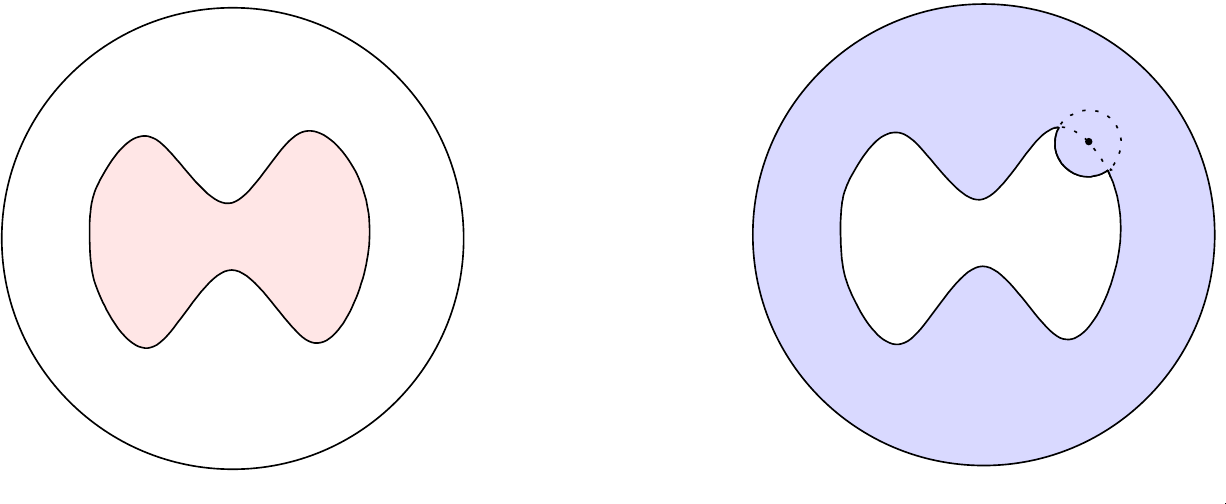}%
\end{picture}%
\setlength{\unitlength}{4144sp}%
\begingroup\makeatletter\ifx\SetFigFont\undefined%
\gdef\SetFigFont#1#2#3#4#5{%
  \reset@font\fontsize{#1}{#2pt}%
  \fontfamily{#3}\fontseries{#4}\fontshape{#5}%
  \selectfont}%
\fi\endgroup%
\begin{picture}(5602,2291)(1,-1472)
\put(4397,324){\makebox(0,0)[lb]{\smash{{\SetFigFont{10}{12.0}{\familydefault}{\mddefault}{\updefault}{\color[rgb]{0,0,0}$V$}%
}}}}
\put(1180,-275){\makebox(0,0)[lb]{\smash{{\SetFigFont{10}{12.0}{\familydefault}{\mddefault}{\updefault}{\color[rgb]{0,0,0}$L_{\text{nt}}$}%
}}}}
\put(643,-283){\makebox(0,0)[lb]{\smash{{\SetFigFont{10}{12.0}{\familydefault}{\mddefault}{\updefault}{\color[rgb]{0,0,0}$U$}%
}}}}
\put(4861, 74){\makebox(0,0)[lb]{\smash{{\SetFigFont{10}{12.0}{\familydefault}{\mddefault}{\updefault}{\color[rgb]{0,0,0}$p$}%
}}}}
\put(5131,119){\makebox(0,0)[lb]{\smash{{\SetFigFont{10}{12.0}{\familydefault}{\mddefault}{\updefault}{\color[rgb]{0,0,0}$U_p$}%
}}}}
\end{picture}%

\end{center}

In dimension $n\geqslant 3$, the gluing lemma provides a way to extend holomorphic line bundles different from the method in Theorem \ref{FSW}.

In \cite{FSW}, the strongly psh exhaustion function is modified to become a nice Morse exhaustion function, also denoted by $\rho$. For any $a\in\br$ and any holomorphic line bundle $L^a$, defined over the super level set $\varOmega^a:=\{\rho>a\}$, they proved that for any point $p$ in the level set $\Gamma^a:=\{\rho=a\}$, there exists a small neighborhood $U_p\subset\varOmega$ of $p$ such that $L^a|_{U_p\cap\varOmega^a}$ is trivial and $L^a$ can be extended trivially to $\varOmega^a\cup U_p$, no matter $p$ is a critical point of the Morse function $\rho$ or not. Since the level set $\Gamma^a$ is compact, after finitely many steps, $L^a$ extends as $L^b$ over $\varOmega^b$ with some $b<a$.

Keep extending $L^a$ until a local minimum $q$ of $\rho$ is reached. The minimum $q$ is an isolated point. There exists some small punctured ball $B_q^*$ centered at $q$ such that $H^1(B_q^*,\co^*)=0$ when $n\geqslant 3$, by a special case of Andreotti-Grauert theory, Proposition 12 in \cite{AG}, which is also proved in \cite{Scheja}. Thus one can extend any holomorphic line bundle trivially across any such local minimum. This proves Theorem \ref{FSW}.

The crucial point above is the following uniqueness result, which a consequence of \cite{FSW}. Let us call it `{\sl downward uniqueness}', since the isomorphism passes to a lower super level set.

\begin{Prop} {\em (Downward uniqueness)} If $L,L'$ are two holomorphic line bundles defined over $\varOmega^b$ that are isomorphic over $\varOmega^a\neq\emptyset$ with $b<a$, then they are isomorphic over $\varOmega^b$.
\end{Prop}

\begin{center}
\begin{picture}(0,0)%
\includegraphics{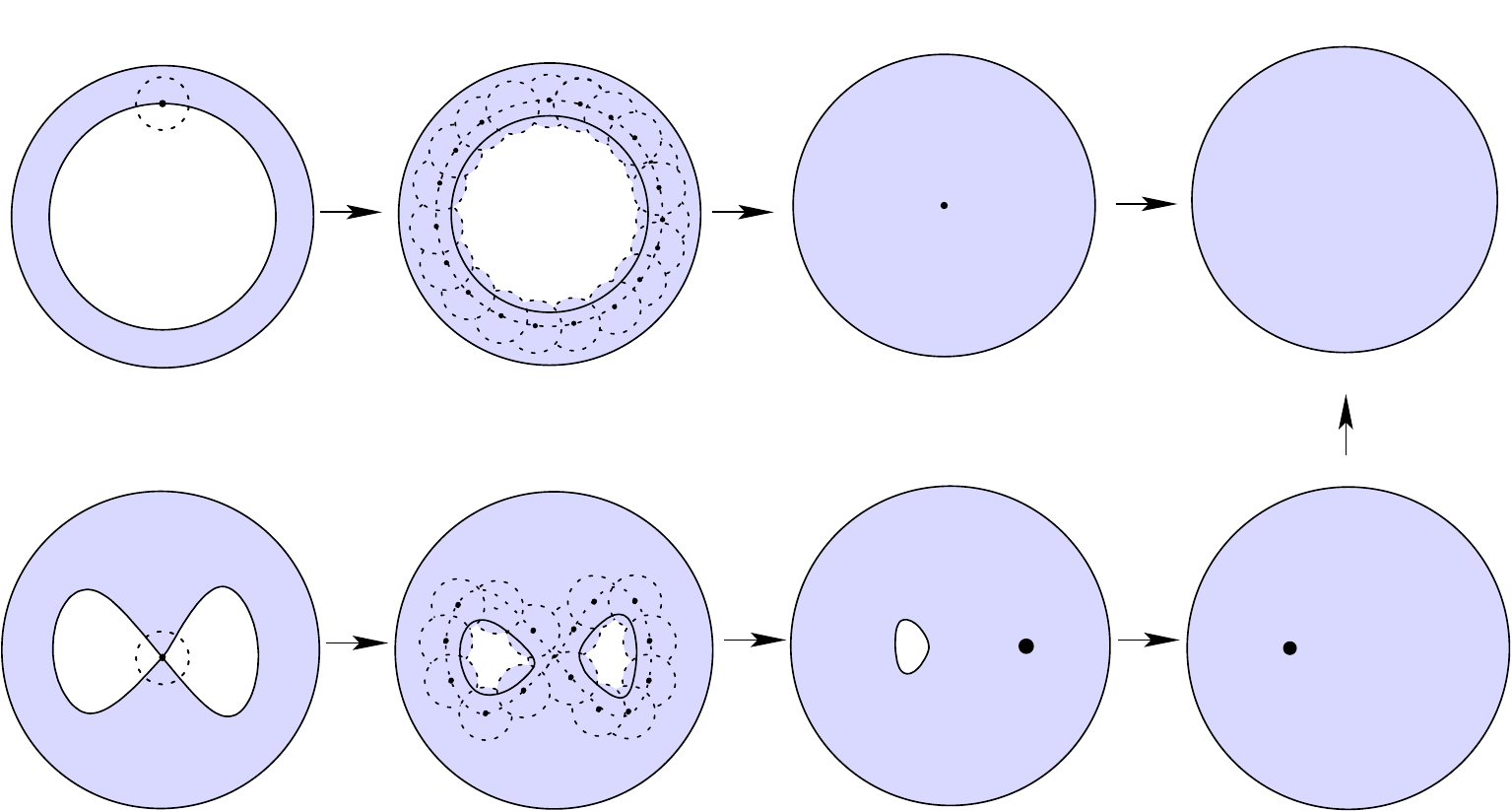}%
\end{picture}%
\setlength{\unitlength}{4144sp}%
\begingroup\makeatletter\ifx\SetFigFont\undefined%
\gdef\SetFigFont#1#2#3#4#5{%
  \reset@font\fontsize{#1}{#2pt}%
  \fontfamily{#3}\fontseries{#4}\fontshape{#5}%
  \selectfont}%
\fi\endgroup%
\begin{picture}(7014,3757)(11,-2898)
\put(2247,-298){\makebox(0,0)[lb]{\smash{{\SetFigFont{10}{12.0}{\familydefault}{\mddefault}{\updefault}{\color[rgb]{0,0,0}$\Gamma^b$}%
}}}}
\put(691,157){\makebox(0,0)[lb]{\smash{{\SetFigFont{10}{12.0}{\familydefault}{\mddefault}{\updefault}{\color[rgb]{0,0,0}$p$}%
}}}}
\put(150,712){\makebox(0,0)[lb]{\smash{{\SetFigFont{10}{12.0}{\familydefault}{\mddefault}{\updefault}{\color[rgb]{0,0,0}$\text{across  a regular point}$}%
}}}}
\put(101,-1354){\makebox(0,0)[lb]{\smash{{\SetFigFont{10}{12.0}{\familydefault}{\mddefault}{\updefault}{\color[rgb]{0,0,0}$\text{across a critical point}$}%
}}}}
\put(678,-2414){\makebox(0,0)[lb]{\smash{{\SetFigFont{10}{12.0}{\familydefault}{\mddefault}{\updefault}{\color[rgb]{0,0,0}$p$}%
}}}}
\put(890,-2203){\makebox(0,0)[lb]{\smash{{\SetFigFont{10}{12.0}{\familydefault}{\mddefault}{\updefault}{\color[rgb]{0,0,0}$U_p$}%
}}}}
\put(4332,  5){\makebox(0,0)[lb]{\smash{{\SetFigFont{10}{12.0}{\familydefault}{\mddefault}{\updefault}{\color[rgb]{0,0,0}$q$}%
}}}}
\put(4688,-2021){\makebox(0,0)[lb]{\smash{{\SetFigFont{10}{12.0}{\familydefault}{\mddefault}{\updefault}{\color[rgb]{0,0,0}$q_1$}%
}}}}
\put(5880,-2028){\makebox(0,0)[lb]{\smash{{\SetFigFont{10}{12.0}{\familydefault}{\mddefault}{\updefault}{\color[rgb]{0,0,0}$q_2$}%
}}}}
\put(4156,-2444){\makebox(0,0)[lb]{\smash{{\SetFigFont{10}{12.0}{\familydefault}{\mddefault}{\updefault}{\color[rgb]{0,0,0}$\Gamma^c$}%
}}}}
\put(2211,-2750){\makebox(0,0)[lb]{\smash{{\SetFigFont{10}{12.0}{\familydefault}{\mddefault}{\updefault}{\color[rgb]{0,0,0}$\Gamma^b$}%
}}}}
\put(310,-2633){\makebox(0,0)[lb]{\smash{{\SetFigFont{10}{12.0}{\familydefault}{\mddefault}{\updefault}{\color[rgb]{0,0,0}$\Gamma^a$}%
}}}}
\put(398,-390){\makebox(0,0)[lb]{\smash{{\SetFigFont{10}{12.0}{\familydefault}{\mddefault}{\updefault}{\color[rgb]{0,0,0}$\Gamma^a$}%
}}}}
\end{picture}%

\end{center}

However, by the gluing Lemma \ref{gluing-lemma}, we can lose uniqueness when we extend through compact sets having shapes different from $K_a$. In our cex constructed in Section \ref{sect-counter-dim-3}, the ball $\varOmega=B(2\sqrt{n}\, e^{\sqrt{\epsilon}})$ admits a strongly psh exhaustion function $\rho(z)=-\log(d(z,\partial\varOmega))$. Since $K$ is compact in $\varOmega$, there exists some $a\in \br$ such that $\varOmega^a=\{\rho>a\}\subset \varOmega\backslash K$. We could restrict the non-extendable holomorphic line bundle $L_{\text{cex}}$, mentioned above, to $\varOmega^a$, and extend $L_{\text{cex}}|_{\varOmega^a}$ to $\varOmega$ by Theorem \ref{FSW}. But in this way we will get a trivial line bundle, which does not agree with the initial bundle $L_{\text{cex}}$ over $\varOmega\backslash K$. In other words, we have the following commutative diagram of restriction maps that are group homomorphisms
\begin{diagram}
H^1(\varOmega,\co^*) & \rTo^{\res^4} & H^1(\varOmega\backslash K,\co^*)\\
& \rdTo^{\res^6} & \dTo^{\res^5}\\
& & H^1(\varOmega^a,\co^*).
\end{diagram}
By Theorem \ref{FSW}, the map $\res^6$ is bijective. But $\res^5$ is not injective since the nontrivial line bundle $L_{\text{cex}}$ and a trivial one $L_{\text{triv}}$ over $\varOmega\backslash K$ have the same restriction on $\varOmega^a$. Consequently $\res^4$ is {\em not} surjective.
\begin{center}
\begin{picture}(0,0)%
\includegraphics{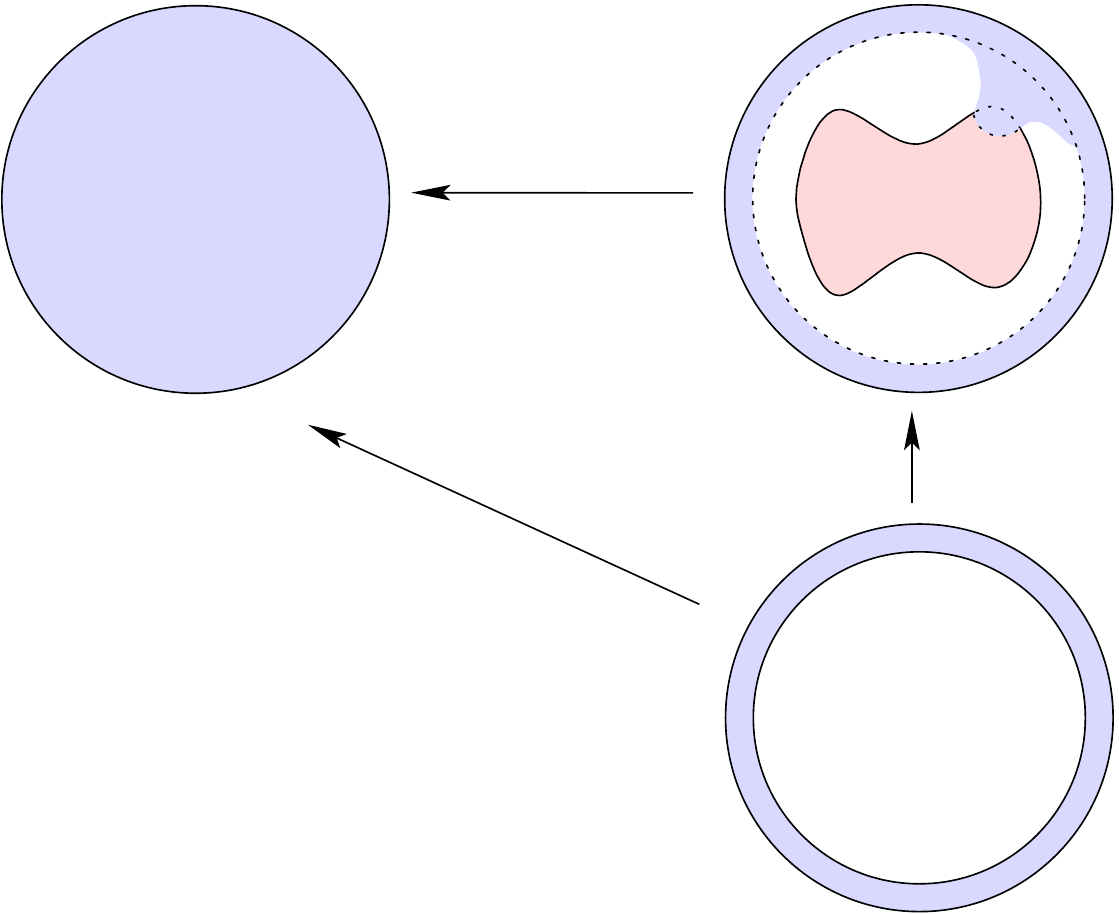}%
\end{picture}%
\setlength{\unitlength}{4144sp}%
\begingroup\makeatletter\ifx\SetFigFont\undefined%
\gdef\SetFigFont#1#2#3#4#5{%
  \reset@font\fontsize{#1}{#2pt}%
  \fontfamily{#3}\fontseries{#4}\fontshape{#5}%
  \selectfont}%
\fi\endgroup%
\begin{picture}(5097,4160)(-5,-3309)
\put(2064,125){\makebox(0,0)[lb]{\smash{{\SetFigFont{10}{12.0}{\familydefault}{\mddefault}{\updefault}{\color[rgb]{0,0,0}$\text{non extendable}$}%
}}}}
\put(4486,413){\makebox(0,0)[lb]{\smash{{\SetFigFont{10}{12.0}{\familydefault}{\mddefault}{\updefault}{\color[rgb]{0,0,0}$L_{\text{triv}}$}%
}}}}
\put(809,-45){\makebox(0,0)[lb]{\smash{{\SetFigFont{10}{12.0}{\familydefault}{\mddefault}{\updefault}{\color[rgb]{0,0,0}$\varOmega$}%
}}}}
\put(4331,-94){\makebox(0,0)[lb]{\smash{{\SetFigFont{10}{12.0}{\familydefault}{\mddefault}{\updefault}{\color[rgb]{0,0,0}$L_{\text{nt}}$}%
}}}}
\put(3811,-90){\makebox(0,0)[lb]{\smash{{\SetFigFont{10}{12.0}{\familydefault}{\mddefault}{\updefault}{\color[rgb]{0,0,0}$U$}%
}}}}
\put(2184,-1310){\makebox(0,0)[lb]{\smash{{\SetFigFont{10}{12.0}{\familydefault}{\mddefault}{\updefault}{\color[rgb]{0,0,0}$\text{FSW extension}$}%
}}}}
\put(4371,-1312){\makebox(0,0)[lb]{\smash{{\SetFigFont{10}{12.0}{\familydefault}{\mddefault}{\updefault}{\color[rgb]{0,0,0}$\text{extension by gluing}$}%
}}}}
\put(4107,-2061){\makebox(0,0)[lb]{\smash{{\SetFigFont{10}{12.0}{\familydefault}{\mddefault}{\updefault}{\color[rgb]{0,0,0}$\varOmega^a$}%
}}}}
\end{picture}%

\end{center}

In conclusion, the map~{\thetag{\ref{H-1-O-star-Omega-K}}}
is {\em not} always surjective, in any dimension $n \geqslant 2$.

\medskip\noindent{\bf Acknowledgments.} The author adresses sincere thanks to Jo\"el Merker for driving him to this problem and for useful discussions. The author also thanks an anonymous referee for pointing out minor mistakes in the first version.

\section{Background}
\label{sec-back}

Now we present the ingredients (1) and (2) mentioned in the Introduction.

\begin{Thm} {\bf[Cartan's theorem B]} Let $X$ be a Stein manifold, $\cf$ be a coherent analytic sheaf on $X$. Then
\[
H^r(X,\cf)=0
\quad{\scriptstyle (r\geqslant 1)}.
\]
\end{Thm}

A proof can be found in Cartan's original paper \cite{Cartan}. Recall the exponential sequence
\[0\rw\bz\xrw{\times 2\pi i}\co\xrw{\exp}\co^*\rw 0\] of sheaves over $X$ induces an exact sequence of cohomologies
\al{H^1(X,\co)\rw H^1(X,\co^*)\rw H^2(X,\bz).}

When $X$ is stein, by Cartan's theorem B we know $H^1(X,\co)=0$. Moreover, if $H^2(X,\bz)=0$, for example when $X$ is contractible, then we get $H^1(X,\co^*)=0$. So we have the following criterion:

\begin{Cor}\label{crit-1} Every holomorphic line bundle over a Stein contractible manifold is trivial.
\end{Cor}

In particular, every convex domain in $\bc^n$ {\small($n\geqslant 1$)} is Stein and contractible.

\begin{Cor}\label{crit-2} Every holomorphic line bundle over a convex domain in $\bc^n$ {\small($n\geqslant 1$)} is trivial.
\end{Cor}

The next ingredient is the extension of holomorphic functions across a totally real plane.

\begin{Thm} Let $D\subset\bc^n$ {\small($n\geqslant 2$)} be a domain, $K\subset \bc^n$ be a totally real plane. Then the restriction map
\al{H^0(D,\co)\rw H^0(D\backslash K,\co)}
is bijective.
\end{Thm}

A proof can be found in the first Chapter of Siu's book \cite{Siu}. We can also apply the argument in the proof of Corollary \ref{Hart-2}.

\begin{Cor}\label{Hart-4} Under the same assumptions, the restriction map
\[
H^0(D,\co^*)\rw H^0(D\backslash K,\co^*)
\]
is bijective.
\end{Cor}

\section{Compactification of Ivashkovich's Counterexample}
\label{sect-counter-dim-2}

In this section we construct a cex in dimension 2. Let $z_1=x_1+iy_1$, $z_2=x_2+iy_2$ be the standard coordinates of $\bc^2\cong\br^4$. For any $r>0$, let $D_r:=\{|x_1|<r,|x_2|<r,y_1^2+y_2^2<1\}$ be a convex bounded domain in $\bc^2$. Let $K:=\{y_1=y_2=0\}$ be a totally real plane in $\bc^2$.  Since $D_r\backslash K$ contracts to $C:=\{x_1=x_2=0,y_1^2+y_2^2=1/2\}\cong S^1$, we have $H^1(D_r\backslash K,\bz)=\bz.$

We can represent a generator of this free $\bz$-module explicitly by using the \v{C}ech cohomology. Take an open cover $\cu=\{U_1,U_2\}$ of $D_r\backslash K$ with $U_1:=\{y_2<|y_1|\}\cap D_r$ and $U_2:=\{y_2>-|y_1|\}\cap D_r$. Then $U_1\cap U_2$ has 2 components, $U_{12}^1:=\{y_1<-|y_2|\}\cap D_r$ and $U_{12}^2:=\{y_1>|y_2|\}\cap D_r$. Let $c=\{c_{12}\}\in Z^1(\cu,\bz)$ be a \v{C}ech 1-cocycle defined by $c_{12}|_{U_{12}^1}=0$, $c_{12}|_{U_{12}^2}=1.$ Then $c$ is a nontrivial 1-cocycle representing a generator $[c]$ of the $\bz-$module $H^1(D_r\backslash K,\bz)$.

\begin{center}
\begin{picture}(0,0)%
\includegraphics{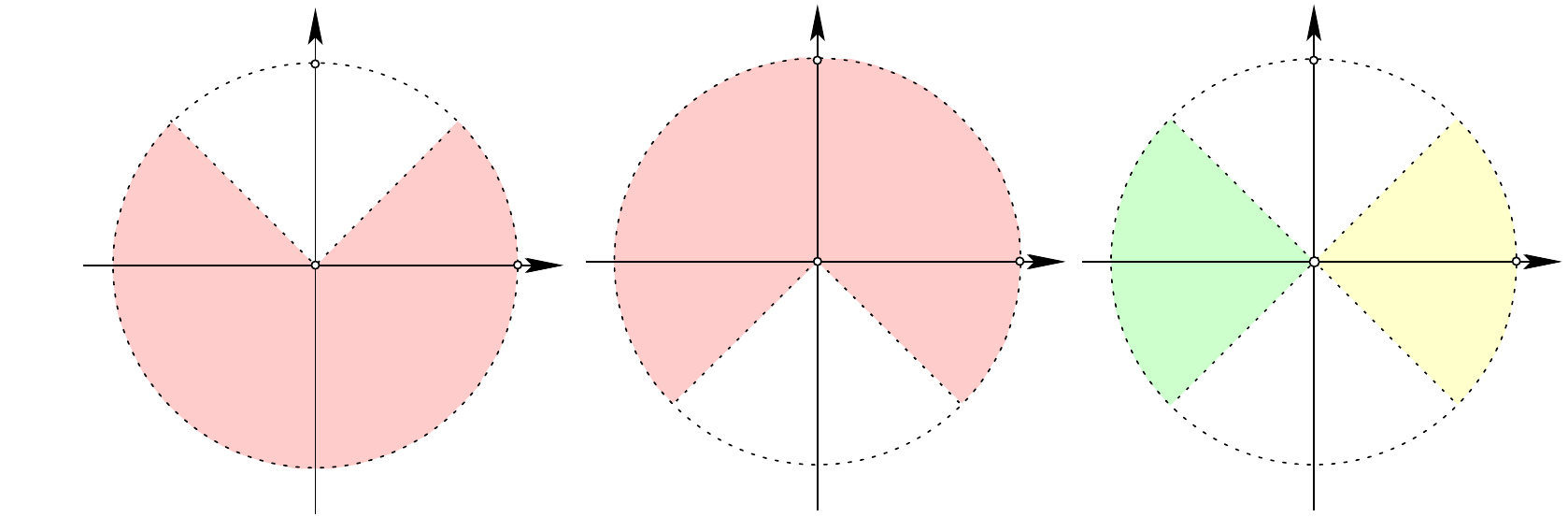}%
\end{picture}%
\setlength{\unitlength}{4144sp}%
\begingroup\makeatletter\ifx\SetFigFont\undefined%
\gdef\SetFigFont#1#2#3#4#5{%
  \reset@font\fontsize{#1}{#2pt}%
  \fontfamily{#3}\fontseries{#4}\fontshape{#5}%
  \selectfont}%
\fi\endgroup%
\begin{picture}(7656,2523)(-372,-1738)
\put(3708,638){\makebox(0,0)[lb]{\smash{{\SetFigFont{10}{12.0}{\familydefault}{\mddefault}{\updefault}{\color[rgb]{0,0,0}$y_2$}%
}}}}
\put(1215,303){\makebox(0,0)[lb]{\smash{{\SetFigFont{10}{12.0}{\familydefault}{\mddefault}{\updefault}{\color[rgb]{0,0,0}1}%
}}}}
\put(1251,620){\makebox(0,0)[lb]{\smash{{\SetFigFont{10}{12.0}{\familydefault}{\mddefault}{\updefault}{\color[rgb]{0,0,0}$y_2$}%
}}}}
\put(7079,-379){\makebox(0,0)[lb]{\smash{{\SetFigFont{10}{12.0}{\familydefault}{\mddefault}{\updefault}{\color[rgb]{0,0,0}$y_1$}%
}}}}
\put(7090,-703){\makebox(0,0)[lb]{\smash{{\SetFigFont{10}{12.0}{\familydefault}{\mddefault}{\updefault}{\color[rgb]{0,0,0}1}%
}}}}
\put(4650,-379){\makebox(0,0)[lb]{\smash{{\SetFigFont{10}{12.0}{\familydefault}{\mddefault}{\updefault}{\color[rgb]{0,0,0}$y_1$}%
}}}}
\put(4661,-703){\makebox(0,0)[lb]{\smash{{\SetFigFont{10}{12.0}{\familydefault}{\mddefault}{\updefault}{\color[rgb]{0,0,0}1}%
}}}}
\put(6101,321){\makebox(0,0)[lb]{\smash{{\SetFigFont{10}{12.0}{\familydefault}{\mddefault}{\updefault}{\color[rgb]{0,0,0}1}%
}}}}
\put(6137,638){\makebox(0,0)[lb]{\smash{{\SetFigFont{10}{12.0}{\familydefault}{\mddefault}{\updefault}{\color[rgb]{0,0,0}$y_2$}%
}}}}
\put(2193,-397){\makebox(0,0)[lb]{\smash{{\SetFigFont{10}{12.0}{\familydefault}{\mddefault}{\updefault}{\color[rgb]{0,0,0}$y_1$}%
}}}}
\put(2204,-721){\makebox(0,0)[lb]{\smash{{\SetFigFont{10}{12.0}{\familydefault}{\mddefault}{\updefault}{\color[rgb]{0,0,0}1}%
}}}}
\put(1228,-709){\makebox(0,0)[lb]{\smash{{\SetFigFont{10}{12.0}{\familydefault}{\mddefault}{\updefault}{\color[rgb]{0,0,0}$U_1$}%
}}}}
\put(3672,321){\makebox(0,0)[lb]{\smash{{\SetFigFont{10}{12.0}{\familydefault}{\mddefault}{\updefault}{\color[rgb]{0,0,0}1}%
}}}}
\put(1921,259){\makebox(0,0)[lb]{\smash{{\SetFigFont{10}{12.0}{\familydefault}{\mddefault}{\updefault}{\color[rgb]{0,0,0}$\{y_2=y_1\}$}%
}}}}
\put(-357,254){\makebox(0,0)[lb]{\smash{{\SetFigFont{10}{12.0}{\familydefault}{\mddefault}{\updefault}{\color[rgb]{0,0,0}$\{y_2=-y_1\}$}%
}}}}
\put(4362,-1338){\makebox(0,0)[lb]{\smash{{\SetFigFont{10}{12.0}{\familydefault}{\mddefault}{\updefault}{\color[rgb]{0,0,0}$\{y_2=-y_1\}$}%
}}}}
\put(2252,-1372){\makebox(0,0)[lb]{\smash{{\SetFigFont{10}{12.0}{\familydefault}{\mddefault}{\updefault}{\color[rgb]{0,0,0}$\{y_2=y_1\}$}%
}}}}
\put(5229,-667){\makebox(0,0)[lb]{\smash{{\SetFigFont{10}{12.0}{\familydefault}{\mddefault}{\updefault}{\color[rgb]{0,0,0}$U_{12}^1$}%
}}}}
\put(6440,-663){\makebox(0,0)[lb]{\smash{{\SetFigFont{10}{12.0}{\familydefault}{\mddefault}{\updefault}{\color[rgb]{0,0,0}$U_{12}^2$}%
}}}}
\put(3670,-400){\makebox(0,0)[lb]{\smash{{\SetFigFont{10}{12.0}{\familydefault}{\mddefault}{\updefault}{\color[rgb]{0,0,0}$U_2$}%
}}}}
\end{picture}%

\end{center}

Recall that the exponential sequence
\al{0\rw\bz\xrw{\times 2\pi i}\co\xrw{\exp}\co^*\rw 0}
induces a long exact sequence
\begin{align}\label{long}
H^0(D_r\backslash K,\co)\xrw{\alpha} H^0(D_r\backslash K,\co^*)\xrw{\beta} H^1(D_r\backslash K,\bz)\xrw{\gamma} H^1(D_r\backslash K,\co)\xrw{\delta} H^1(D_r\backslash K,\co^*).
\end{align}

By Corollary \ref{Hart-4}
\[
H^0(D_r\backslash K,\co^*)\cong H^0(D_r,\co^*).
\]
Since $D_r$ is simply connected, the map
\al{H^0(D_r,\co)\xrw{\exp}H^0(D_r,\co^*)}
is surjective, thus in the long exact sequence (\ref{long}), the map
\al{\alpha:H^0(D_r\backslash K,\co)\rw H^0(D_r\backslash K,\co^*)}
is also surjective. We have $\im(\alpha)=\ker(\beta)=H^0(D_r\backslash K,\co^*)$ so $0=\im(\beta)=\ker(\gamma)$, i.e. $\gamma$ is injective. We know that the sequence
\al{0\rw H^1(D_r\backslash K,\bz)\xrw{\gamma} H^1(D_r\backslash K,\co)\xrw{\delta} H^1(D_r\backslash K,\co^*)}
is exact.

The generator $[c]\in H^1(D_r\backslash K,\bz)$ maps to $\gamma([c])$ which is nontrivial since $\gamma$ is injective. Since $\frac{1}{2}\gamma([c])$ is not in the image of $\gamma$ we know $\delta\big(\frac{1}{2} \gamma([c])\big)$ represents a nontrivial holomorphic line bundle over $D_r\backslash K$, which cannot be extended to $D_r$ since $H^1(D_r,\co^*)=0$.

Using \v{C}ech cohomology, the element $\gamma([c])$ can be represented by $\{2\pi i c_{12}\}\in Z^1(\cu,\co)$ and $\delta\big(\frac{1}{2}\gamma([c])\big)$ can be represented by $\{e^{\pi ic_{12}}\}\in Z^1(\cu,\co^*)$. Denote this 1-cocycle by $f=\{f_{12}\}$. We have $f_{12}|_{U_{12}^1}=1$, $f_{12}|_{U_{12}^2}=-1$. Let $L$ be a holomorphic line bundle defined on $D_r\backslash K$, trivial on $U_1$ and $U_2$ and the transition function is defined by $\{f_{12}\}$. Then $L$ is a nontrivial holomorphic line bundle. $L$ cannot be extended to $D_r$, since by Corollary \ref{crit-2} every holomorphic line bundle over $D_r$ is trivial.

\medskip

Now we will construct the following objects:
\begin{itemize}
\item a bounded pseudoconvex domain $\varOmega\subset\bc^2$ with a strongly $\mathscr{C}^\infty$ psh exhaustion function $\rho$;
\item some $a\in\br$ and some compact $K_a:=\rho^{-1}(-\infty,a]\subset\subset \varOmega$;
\item a holomorphic line bundle $L$ over $\varOmega\backslash K_a$ which can not be extended to $\varOmega$.
\end{itemize}

Recall $D_r=\{|x_1|,|x_2|<r,y_1^2+y_2^2<1\}$. Consider the map
\al{\varphi=\exp(i\cdot):\bc^2&\rw \bc^2,\\
(z_1,z_2)&\longmapsto (w_1,w_2):=(e^{iz_1},e^{iz_2}).}
This map $\varphi$ is locally biholomorphic. It is bijective, hence biholomorphic from $D_r$ onto $\varphi(D_r)$, when $r\leqslant\pi$. When $r>\pi$, the image is $\{(\log|w_1|)^2+(\log|w_2|)^2<1\}$. We define $\varOmega$ as this open set. Actually $\varOmega$ is a Grauert tube around the totally real torus $\{|w_1|=|w_2|=1\}$. It is a bounded pseudoconvex domain, as a special case of Proposition \ref{prop:G} with $n=2$ and $\epsilon=1$. The function
\al{\rho:\varOmega&\longrightarrow[0,+\infty)\\
(w_1,w_2)&\longmapsto-\log\big(1-(\log|w_1|)^2-(\log|w_2|)^2\big)}
is a strongly psh exhaustion function of $\varOmega$ and $\varphi(K\cap D_r)=\{|w_1|=|w_2|=1\}=\rho^{-1}(0)=K_0\subset\subset \varOmega$.

Recall the covering $\cu$ of $D_r\backslash K$ and the \v{C}ech 1-cocycle $\{f_{12}\}$ above. Notice that $f_{12}$ is constant along the $(x_1,x_2)$-directions. In particular,
\al{f_{12}(x_1+2k_1\pi,x_2+2k_2\pi,y_1,y_2)=f_{12}(x_1,x_2,y_1,y_2)}
whenever $k_1$, $k_2\in\bz$, $(x_1,x_2,y_1,y_2)\in D_r$ and $(x_1+2k_1\pi,x_2+2k_2\pi,y_1,y_2)\in D_r$. So $f_{12}$ induces a function $\tilde{f}_{12}$ well defined on the disjoint union $\varphi(U_{12})=\varphi(U_{12}^1)\sqcup\varphi(U_{12}^2)$ with 
\al{\tilde{f}_{12}|_{\varphi(U_{12}^1)}&=1,\\
\tilde{f}_{12}|_{\varphi(U_{12}^2)}&=-1.}
Here $\varphi(\cu):=\{\varphi(U_1),\varphi(U_2)\}$ is an open cover of $\varOmega\backslash K_0$. This open cover, together with the transition function $\tilde{f}_{12}$, defines a nontrivial holomorphic line bundle $L$ over $\varOmega\backslash K_0$.

Suppose $L$ can be extended to a holomorphic line bundle $L'$ over $\varOmega$. Note that $\varphi=\exp(i\cdot)$ is a biholomorphism between $D_1=\{|x_1|,|x_2|<1,y_1^2+y_2^2<1\}$ and its image $\varphi(D_1)$. The pull-back $\varphi^*(L'|_{\varphi(D_1)})$ gives a holomorphic line bundle defined over $D_1$ which extends $\varphi^*(L|_{\varphi(D_1)\backslash K_0})$. Here $\varphi^*(L|_{\varphi(D_1)\backslash K_0})$ is a holomorphic line bundle defined over $D_1\backslash K$, since $\varphi^{-1}\big(\varphi(D_1)\backslash K_0\big)=D_1\backslash K$. However, due to our discussion in Section 2, such $\varphi^*(L|_{\varphi(D_1)\backslash K_0})$ is nontrivial hence cannot be extended across $K$. This contradiction shows that $L$ cannot be extended to $\varOmega$.

We can draw $\varOmega\subset \br^4$ as a movie of its 3d-sections (when $y_2$ is fixed) in $\br^3$.

\begin{figure}[h]
    \centering 
    \includegraphics[width=0.7\textwidth]{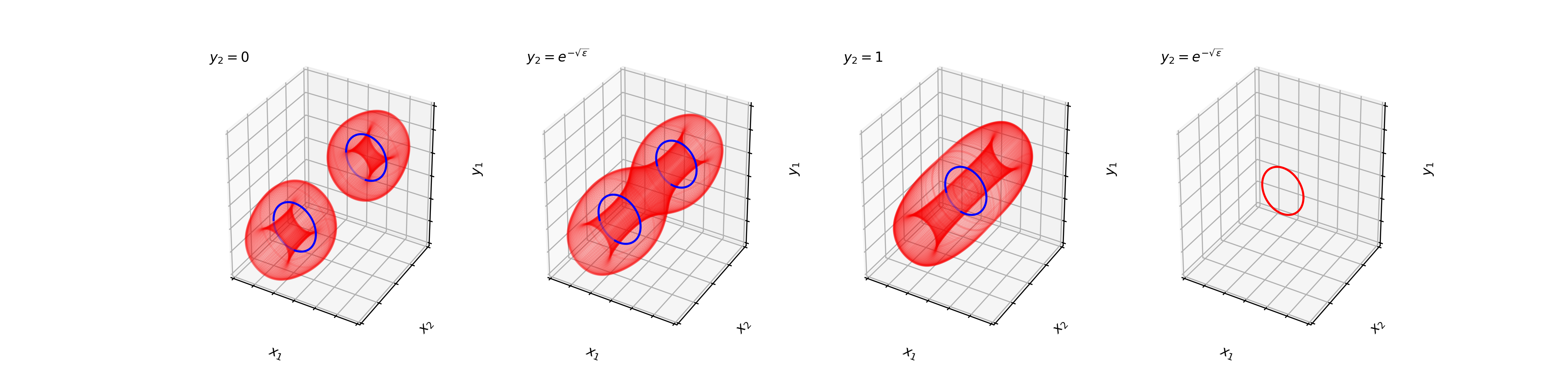} 
  \caption{3d-sections of $\varOmega$, where the red surface is $\partial\varOmega$ and the blue curve is $K_0$} 
\end{figure}

In fact, each 3d-section is obtained by rotating the 2d-section (when $y_1=0$) along the $x_2$-axis (the dashed line).

\begin{figure}[htbp]
    \centering 
    \includegraphics[width=0.8\textwidth]{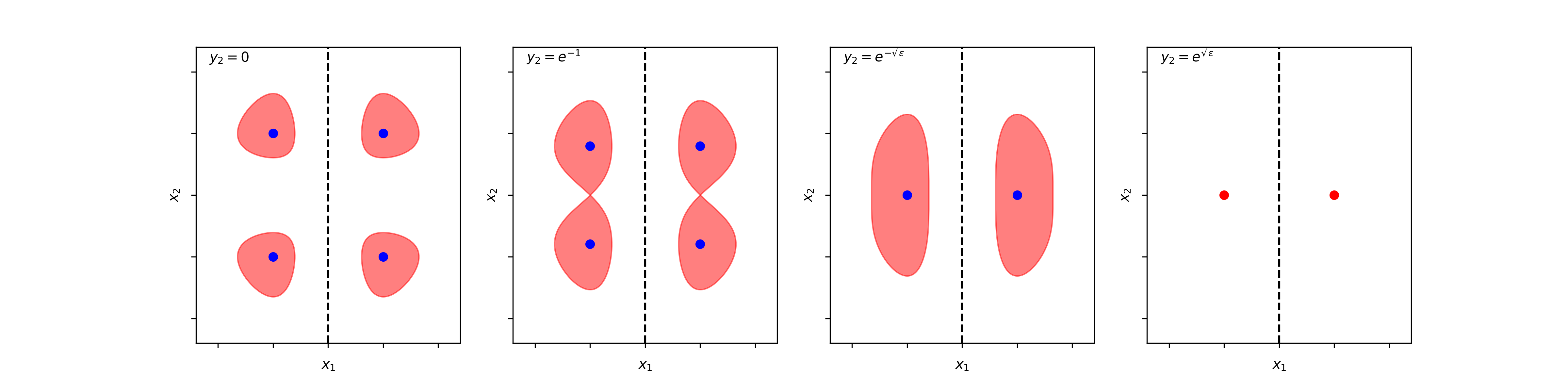} 
  \caption{2d-sections of $\varOmega$ when $y_1=0$} 
\end{figure}

We can also draw $\varphi(U_1),\varphi(U_2),\varphi(U_{12}^1)$ and $\varphi(U_{12}^2)$ in this way.

\begin{figure}[htbp]
    \centering 
    \includegraphics[width=0.8\textwidth]{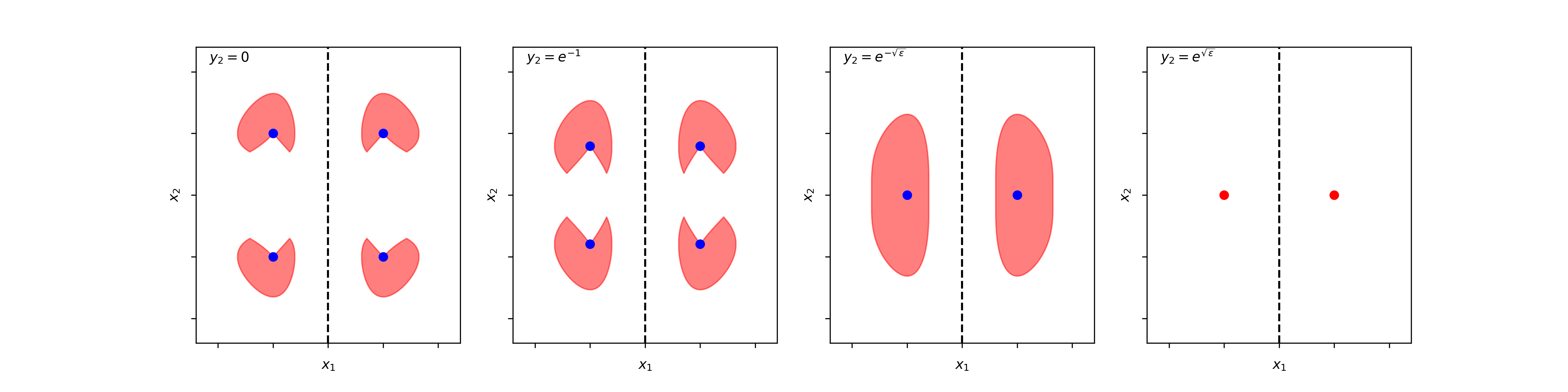} 
  \caption{2d-sections of $\varphi(U_1)$} 
\end{figure}

\begin{figure}[htbp]
    \centering 
    \includegraphics[width=0.8\textwidth]{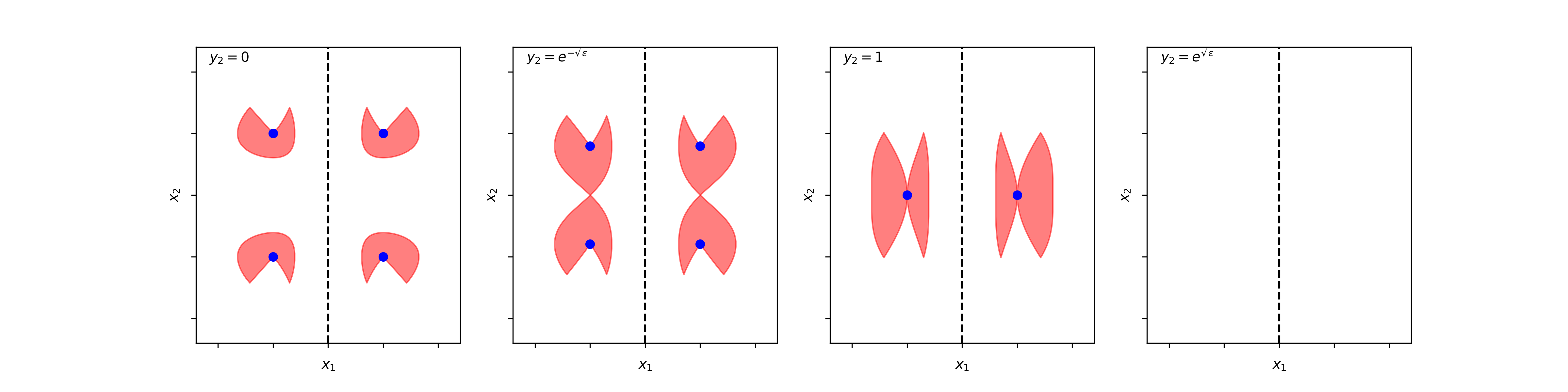} 
  \caption{2d-sections of $\varphi(U_2)$} 
\end{figure}

\begin{figure}[htbp]
    \centering 
    \includegraphics[width=0.8\textwidth]{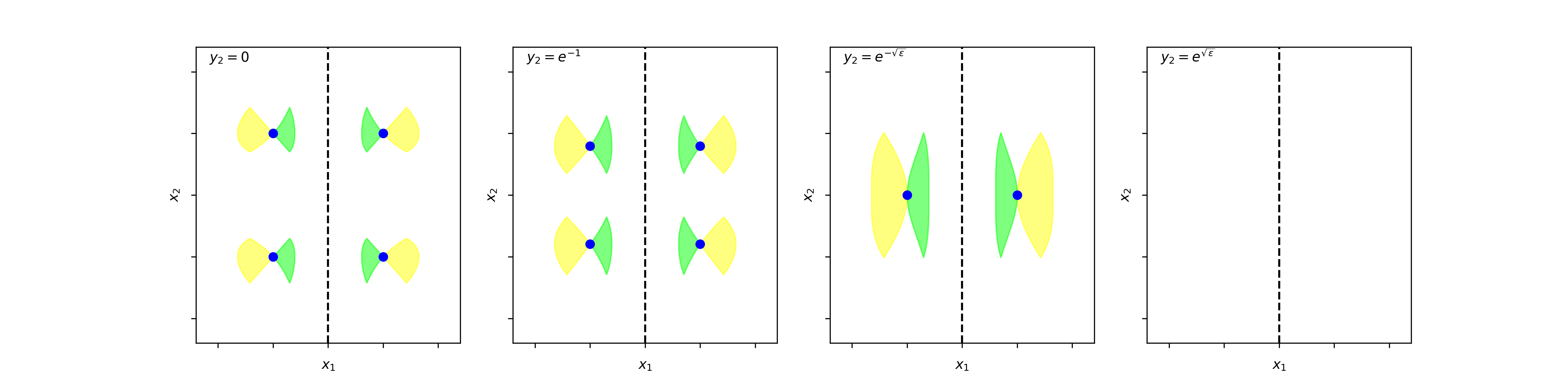} 
  \caption{2d-sections of $\varphi(U_{12}^1)$ in green and of $\varphi(U_{12}^2)$ in yellow} 
\end{figure}

\clearpage
\section{Counterexamples in general dimension}
\label{sect-counter-dim-3}

As announced in the Introduction, we will construct some `strange' bundles which cannot be extended from $\varOmega\backslash K$ to $\varOmega$. The key idea is a certain gluing lemma describing `flexibility' of holomorphic line bundles. Actually such lemma holds for holomorphic vector bundles.

\subsection{Gluing lemma}\label{subsection-gluing}
Roughly speaking, holomorphic vector bundles have more `flexibility' than holomorphic functions. Let $\varOmega\subset\bc^n$ {\small($n\geqslant 2$)} be a domain and $U\subset \varOmega$ be a non-empty open subset. If two holomorphic functions $f$ and $g$ in $\varOmega$ are equal over $U$, then they are equal over $\varOmega$. However, if two holomorphic vector bundles $E$ and $F$ over $\varOmega$ are isomorphic over $U$, then they may not be isomorphic over $\varOmega$ in general. For example, let $E$ and $F$ be two non-isomorphic holomorphic line bundles over $\varOmega$. For any $x\in\varOmega$, there exists a neighborhood $U_E$ (resp. $U_F$) of $x$ in $\varOmega$ where $E$ (resp. $F$) is trivial. So $E$ and $F$ are trivial over $U:=U_E\cap U_F$, another neighborhood of $x$ in $\varOmega$. So they are isomorphic over a non-empty open subset $U$ of $\varOmega$.

\begin{Lem}\label{lem:gluing} {\em (Gluing lemma for holomorphic vector bundles)} Let $X$ be a complex manifold, let $U,V\subset X$ be two open subsets and let $W:=U\cap V$. For any integer $r\geqslant 1$, let $E_U$ be a holomorphic vector bundle of rank $r$ over $U$. Let $E_V$ be a holomorphic vector bundle of rank $r$ over $V$ such that $E_U|_W\cong E_V|_W$. Then there exists a holomorphic vector bundle $E$ of rank $r$ over $U\cup V$ such that $E|_U\cong E_U$ and $E|_V\cong E_V$. 
\end{Lem}

\begin{Rmk} \em{Denote by $\fglr$ the sheaf of invertible $r\times r$ matrices with coefficients in the sheaf $\co$ of holomorphic functions. In particular, $\mathcal{GL}_1(\co)=\co^*$. Using the language of category theory, by the universal property of the fibre product, the following commutative diagram 
\begin{diagram}
H^1\big(U\cup V,\fglr\big) & \rTo^{\text{res}} & H^1\big(U,\fglr\big)\\
\dTo^{\text{res}} & \comm& \dTo^{\text{res}}\\
H^1\big(V,\fglr\big) & \rTo^{\text{res}} & H^1\big(U\cap V,\fglr\big)
\end{diagram}
induces a canonical map
\al{q:H^1\big(U\cup V,\fglr\big)\cong H^1\big(U,\fglr\big)\times_{H^1(U\cap V,\fglr)}H^1\big(V,\fglr\big).}
The gluing Lemma \ref{lem:gluing} states that $q$ is an epimorphism, by constructing a left inverse of $q$.}
\end{Rmk}

\begin{Rmk} \em{Note that we only have existence, but not uniqueness in general. That is to say, $E$ is not uniquely determined up to isomorphism by the information of $E_U$ and $E_V$. When $r=1$, a simple cex to uniqueness can be constructed by taking
\al{
&X=\bc\bp^1:=\{[Z_0:Z_1],(Z_0,Z_1)\in\bc^2\backslash(0,0)\},\\
&U:=\{Z_0\neq 0\}=\{[1:{\textstyle{\frac{Z_1}{Z_0}}}],{\textstyle{\frac{Z_1}{Z_0}}}\in\bc\}\cong \bc,\\
&V:=\{Z_1\neq 0\}=\{[{\textstyle{\frac{Z_0}{Z_1}}}:1],{\textstyle{\frac{Z_0}{Z_1}}}\in\bc\}\cong\bc,\\
&W=U\cap V\cong \bc^*
}
and $E_U:=U\times\bc$ with coordinates $(z_1,s_1)$, $E_V:=V\times \bc$ with coordinates $(z_0,s_0)$ being trivial line bundles with the identifications
\al{
(z_1,s_1)=({\textstyle{\frac{1}{z_0}}},{\textstyle{\frac{1}{z_0^n}}}s_0)
}
for each $n\in\bz$. This defines the holomorphic line bundle $O(n)$, trivial over $U$ and $V$. But $O(n)$ and $O(m)$ are not isomorphic whenever $n\neq m$.}
\end{Rmk}

\proof[Proof of Lemma~\ref{lem:gluing}] It suffices to consider the case where $U$, $V$, $W$ are non-empty. Denote the projection maps by $\pi_U:E_U\rw U$ and $\pi_V:E_V\rw V$. Let $\{(U_i,\varphi_{i})\}$ be a trivialization of $E_U$, $\{(V^j,\psi^j)\}$ be a trivialization of $E_V$. We will use the notations
\al{U_{i_1,\dots,i_s}:=\bigcap\limits_{x=1}^s U_{i_x},\\
V^{j_1,\dots,j_t}:=\bigcap\limits_{y=1}^t V^{j_y}.}
By trivializations of these vector bundles, we mean that $\{U_i\}$ is an open cover of $U$, $\{V^j\}$ is an open cover of $V$ and
\al{\varphi_{i}:\pi_U^{-1}(U_i)\rw U_i\times \bc^r,\\
\psi^j:\pi_V^{-1}(V^j)\rw V^j\times \bc^r,}
are homeomorphisms and the transition functions $f_{i_2,i_1}$, $g^{j_2,j_1}$ defined by
\al{\varphi_{i_2}\circ(\varphi_{i_1})^{-1}:U_{i_1,i_2}\times \bc^r&\rw U_{i_1,i_2}\times \bc^r\\
(z,l)&\longmapsto \big(z,f_{i_2,i_1}(z)\cdot l\big),\\
\psi^{j_2}\circ(\psi^{j_1})^{-1}:V^{j_1,j_2}\times \bc^r&\rw V^{j_1,j_2}\times \bc^r\\
(z,l)&\longmapsto \big(z,g^{j_2,j_1}(z)\cdot l\big),}
are holomorphic maps valued in $\fglr(U_{i_1,i_2})$, $\fglr(V^{j_1,j_2})$, satisfying the {\sl cocycle conditions}:
\begin{align}\label{cocy1}
f_{i_3,i_2}(z)\cdot f_{i_2,i_1} (z)=f_{i_3,i_1}(z)\quad {\scriptstyle(z\in U_{i_1,i_2,i_3})},\\
\label{cocy2}
g^{j_3,j_2}(z)\cdot g^{j_2,j_1}(z)=g^{j_3,j_1}(z)\quad {\scriptstyle(z\in V^{j_1,j_2,j_3})}.
\end{align}

In fact, we can use the notation $H^{0,0}_{\bar{\partial}}(U_{I},E_U)$ for the set of holomorphic sections of $E_U$ over $U_{I}$ where $I=\{i_1,\dots,i_s\}$. It is indeed a free $H^0(U_I,\co)$-module of rank $r$. If we use $e_1,\dots,e_r$ as the standard basis of the $\bc$-vector space $\bc^r$, then $\{(\varphi_{i})^{-1}(z,e_1),$ $\dots$ $,(\varphi_{i})^{-1}(z,e_r)\}$ is a set of nowhere vanishing holomorphic sections of $E_U|_{U_i}$ which generate $H^{0,0}_{\bar{\partial}}(U_i,E_U)$. So actually $f_{i_2,i_1}$ is a $H^0(U_{i_1,i_2},\co)$-coefficients invertible linear map such that for $\alpha=1,\dots,r$,
\al{(\varphi_{i_2})^{-1}(z,e_\alpha)=\sum\limits_{\beta=1}^r\big(f_{i_2,i_1}(z)\big)_{\alpha,\beta}(\varphi_{i_1})^{-1}(z,e_\beta)\quad {\scriptstyle(\beta=1,\dots,r;~z\in U_{i_1,i_2})},}
and the cocycle conditions are automatically satisfied.

In the language of \v{C}ech cohomology we would say that $E_U$ is represented by the open cover $\{U_i\}$ of $U$ and the 1-cocycle $\{f_{i_2,i_1}\}\in Z^1\big(\{U_i\},\fglr\big)$. We have similar statements for $\big(E_V,\{V^j\},\{g^{j_2,j_1}\}\big)$.

Then $\{W_i^j:=U_i\cap V^j\}$ is an open cover of $W$ trivializing $E_U$ and $E_V$ simultaneously, i.e. $H^{0,0}_{\bar{\partial}}(W_i^j,E_U)$ \big(resp. $H^{0,0}_{\bar{\partial}}(W_i^j,E_V)$\big) is a rank $r$ free $H^0(W_i^j,\co)$-module generated by $\{(\varphi_{i})^{-1}(z,e_1),$ $\dots$ $,(\varphi_{i})^{-1}(z,e_r)\}$ \big(resp. $\{(\psi^{j})^{-1}(z,e_1),$ $\dots$ $,(\psi^{j})^{-1}(z,e_r)\}$\big). The isomorphism $h:E_U|_W\cong E_V|_W$ induces an isomorphism between the rank $r$ free $H^0(W_i^j,\co)$-modules $H^{0,0}_{\bar{\partial}}(W_i^j,E_U)$ and $H^{0,0}_{\bar{\partial}}(W_i^j,E_V)$. It is determined by $h^j_{~i}\in H^0\big(W_i^j,\fglr\big)$ such that for $\alpha=1,\dots,r$,
\al{(\psi^{j})^{-1}(z,e_\alpha)=\sum\limits_{\beta=1}^r\big(h^j_{~i}(z)\big)_{\alpha,\beta}(\varphi_{i})^{-1}(z,e_\beta)\quad{\scriptstyle(\beta=1,\dots,r;~ z\in W^j_i)}.
}

Use the notation $W_{i_1,\dots, i_s}^{j_1,\dots, j_t}:=(\bigcap\limits_{x=1}^sU_{i_x})\cap(\bigcap\limits_{y=1}^t V^{j_y})$. For any indices $i_1$, $i_2$, $j_1$, $j_2$ we get the {\sl transition equations}
\begin{align}\label{liso}
f_{i_2,i_1}(z)=\big(h_{~~i_2}^{j_2}(z)\big)^{-1}\cdot g^{j_2,j_1}(z)\cdot h_{~~i_1}^{j_1}(z)\quad{\scriptstyle(z\in W_{i_1,i_2}^{j_1,j_2})}.
\end{align}

Now we define $E$, a holomorphic vector bundle over $U\cup V$. Note that $\{U_i,V^j\}$ is actually an open cover of  $U\cup V$. We let $E$ to be trivial on each $U_i$ and $V^j$ and define transition functions $l\in Z^1\big(\{U_i,V^j\},\fglr\big)$ by
\al{l_{i_2,i_1} &:=f_{i_2,i_1},\\
l^{j_2,j_1} &:=g^{j_2,j_1},\\
l_{~i}^j &:=h_{~i}^j.}
The cocycle conditions among ($U_{i_1}$, $U_{i_2}$, $U_{i_3}$) and ($V^{j_1}$, $V^{j_2}$, $V^{j_3}$) are satisfied because of (\ref{cocy1}) and (\ref{cocy2}). We only need to check the cocycle conditions among ($U_{i_1}$, $U_{i_2}$, $V^j$) and ($U_i$, $V^{j_1}$, $V^{j_2}$), i.e.
\al{h_{~i_2}^j(z)\cdot f_{i_2,i_1}(z) &=h_{~i_1}^j(z)\quad{\scriptstyle(z\in W_{i_1,i_2}^j)},\\
g^{j_2,j_1}(z)\cdot h_{~~i}^{j_1}(z) &=h_{~~i}^{j_2}(z)\quad{\scriptstyle(z\in W_i^{j_1,j_2})}.}
But this can be achieved by taking $j_1=j_2=j$ or $i_1=i_2=i$ in (\ref{liso}) and using $f_{i,i}=id$, $g^{j,j}=id$. We have $E|_U\cong E_U$, because both bundles are given by the same transition functions $\{l_{i_2,i_1}=f_{i_2,i_1}\}$ with respect to the same open covering $\{U_i\}$ of $U$. For the same reason $E|_V\cong E_V$. \endproof

\subsection{A Stein manifold with a nontrivial holomorphic line bundle}

Take $\rho:(\bc^*)^n\rw[0,+\infty)$, $(z_1,\dots,z_n)\mapsto\sum\limits_{j=1}^n(\log|z_j|)^2$. For any $\epsilon>0$, define $G_\epsilon:=\rho^{-1}[0,\epsilon)$. Denote by $B(r):=\{\sum\limits_{j=1}^n|z_j|^2<r^2\}$ the open ball centered at the origin of radius $r>0$ in $\bc^n$.
\begin{Prop}\label{prop:G} For every $\epsilon>0$,
\begin{enumerate}[{\bf (i)}]
\item $G_\epsilon\subset B(\sqrt{n}e^{\sqrt{\epsilon}})$ is bounded;
\item $G_\epsilon$ is pseudoconvex with smooth boundary;
\item $G_\epsilon$ is connected. In fact $G_\epsilon$ contracts to a $n$-dimensional torus. Its Picard group is $H^1(G_\epsilon,\co^*)\cong\bz^{\binom{n}{2}}$. In particular, $G_\epsilon$ carries a nontrivial holomorphic line bundle.
\end{enumerate}
\end{Prop}
\proof {\bf (i)} Since
\begin{align*}
  \sum\limits_{j=1}^n(\log|z_j|)^2 <\epsilon
  ~\Longrightarrow~ &(\log|z_j|)^2 <\epsilon\cond{j=1,\dots,n}\\
  \Longrightarrow~ &e^{-\sqrt{\epsilon}}<|z_j| <e^{\sqrt{\epsilon}}\cond{j=1,\dots,n}\\
  \Longrightarrow~ &\sum\limits_{j=1}^n|z_j|^2 <ne^{2\sqrt{\epsilon}},
\end{align*}
we see that $G_\epsilon\subset B(\sqrt{n}e^{\sqrt{\epsilon}})$.

{\bf (ii)} Actually the boundary $\partial G_\epsilon=\rho^{-1}(\epsilon)$. We know that $\rho$ is smooth on $(\bc^*)^n$ and $d\rho=\sum\limits_{j=1}^n\log|z_j|({\textstyle{\frac{dz_j}{z_j}}}+{\textstyle{\frac{d\bar{z}_j}{\bar{z}_j}}})$.
\al{&d\rho=0\\
\Longleftrightarrow~ &|z_j| =1\cond{j=1,\dots,n}\\
\Longleftrightarrow~ &\rho=0.}
So $\epsilon>0$ is a regular value of $\rho$, hence $\partial G_\epsilon$ is smooth.

To prove that $G_\epsilon$ is pseudoconvex, we check the Levi-condition. For any $z\in\partial G_\epsilon$ and $(w_1,\dots,w_n)\in T_z(\bc^*)^n=T_z\bc^n$ satisfying $\sum\limits_{j=1}^n\frac{\partial\rho(z)}{\partial z_j}w_j=0$, we have
\al{\sum\limits_{j,k=1}^n\frac{\partial^2\rho(z)}{\partial z_j\bar{\partial}z_k}w_j\bar{w}_k =\sum\limits_{j=1}^n\frac{|w_j|^2}{2|z_j|^2}\geqslant \frac{1}{2e^{2\sqrt{\epsilon}}}\sum\limits_{j=1}^n|w_j|^2.}
Thus $G_\epsilon$ is pseudoconvex.

{\bf (iii)} We claim that $G_\epsilon$ contracts to $\rho^{-1}(0)=\{|z_j|=1,j=1,\dots,n\}\cong \underbrace{S^1\times \ldots \times S^1}_{n \text{ times}}=T^n$. To show this, we construct a contraction map
\al{H:G_\epsilon\times [0,1] &\rw \bc^n,\\
(z,t) &\longmapsto (|z_1|^{\lambda(z,t)}z_1,\dots,|z_n|^{\lambda(z,t)}z_n),}
where $\lambda(z,t):G_\epsilon\times[0,1]\rw \br$ is a continuous map. Note that
\al{\rho\big(H(z,t)\big)=\big(\lambda(z,t)+1\big)^2\rho(z).}

Then $H$ is a contraction map if for any $z\in G_\epsilon$ and any $t\in[0,1]$ we have
\begin{equation}
\left\{
\begin{aligned}
H(z,0) =z &\Longrightarrow \lambda(z,0)=0,\\
\im(H) \in G_\epsilon &\Longrightarrow \big(\lambda(z,t)+1\big)^2\rho(z)<\epsilon.
\end{aligned}
\right.
\end{equation}

If we take $\lambda(z,t)=-t$, then these conditions are satisfied and $\im\big(H(z,1)\big)=\rho^{-1}(0)=\{|z_j|=1,j=1,\dots,n\}\cong T^n$. We see $G_\epsilon$ is connected since $T^n$ is.

Now we calculate the Picard group $H^1(G_\epsilon,\co^*)$. By {\bf(ii)}, $G_\epsilon$ is Stein and by Cartan's theorem B, we have $H^1(G_\epsilon,\co)=H^2(G_\epsilon,\co)=0$.
Recall that the exponential exact sequence
\al{0\rw \bz\xrw{\times 2\pi i}\co\xrw{exp}\co^*\rw 0}
induces a long exact sequence
\al{H^1(G_\epsilon,\co)\rw H^1(G_\epsilon,\co^*)\rw H^2(G_\epsilon,\bz)\rw H^2(G_\epsilon,\co).}

Since the first and the last term vanish, we have
\al{H^1(G_\epsilon,\co^*)\cong H^2(G_\epsilon,\bz)\cong H^2(T^n,\bz).}

Recall that by using a Mayer-Vietoris sequence we have $H^k(S^1\times M,\bz)\cong H^{k-1}(M,\bz)\oplus H^k(M,\bz)$ for any simplicial complex $M$ and any $k\in\bz_+$. Hence $H^k(T^n,\bz)\cong\bz^{\binom{n}{k}}$. In particular $H^1(G_\epsilon,\co^*)\cong H^2(T^n,\bz)\cong\bz^{\binom{n}{2}}$ is nontrivial since $n\geqslant 2$. \endproof

\subsection{Gluing process}

Now we take $0<\epsilon<n$, for example $\epsilon=n/2$. The boundary $\partial G_\epsilon$ is then given by the equation
\[
\rho(z)=\sum\limits_{j=1}^n(\log|z_j|)^2=\epsilon.
\]

\begin{Prop}\label{real-hessian}
The real Hessian $H\rho$ is positive definite at the point $p=(e^{\sqrt{\epsilon/n}},\dots,e^{\sqrt{\epsilon/n}})\in\partial G_\epsilon$.
\end{Prop}

\proof For $j=1,\dots,n$, we have
\al{
\frac{\partial^2 \rho}{\partial x_j^2} &=\frac{2}{|z_j|^4}\big((y_j^2-x_j^2)\log|z_j|+x_j^2\big),\\
\frac{\partial^2 \rho}{\partial x_j\partial y_j} &=\frac{2x_jy_j}{|z_j|^4}(1-2\log|z_j|),\\
\frac{\partial^2 \rho}{\partial y_j^2} &=\frac{2}{|z_j|^4}\big((x_j^2-y_j^2)\log|z_j|+y_j^2\big).
}

So the real Hessian is
\[
H\rho=
\begin{bmatrix}
    \frac{2}{|z_1|^4}H^{2\times 2}_1 & 0  & \dots  & 0 \\
    0 & \frac{2}{|z_2|^4}H^{2\times 2}_2 &  \dots  & 0 \\
    \vdots & \vdots & \ddots & \vdots \\
    0 & 0 & \dots  & \frac{2}{|z_n|^4}H^{2\times 2}_n
\end{bmatrix},
\]
where
\[
H^{2\times 2}_j=
\begin{bmatrix}
    (y_j^2-x_j^2)\log|z_j|+x_j^2 & x_jy_j(1-2\log|z_j|) \\
    x_jy_j(1-2\log|z_j|) & (x_j^2-y_j^2)\log|z_j|+y_j^2 
\end{bmatrix}\cond{j=1,\dots,n}.
\]

The real Hessian $H\rho$ is positive definite if and only if for all $j$, $H^{2\times 2}_j$ is positive definite. That is equivalent to
\begin{equation*}
\left\{
\begin{aligned}
\text{tr} (H^{2\times 2}_j)>0,\\
\det(H^{2\times 2}_j)>0,
\end{aligned}
\right.
\end{equation*}
for all $j$. The first inequality is achieved since $\text{tr}(H^{2\times 2}_j)=x_j^2+y_j^2>0$ for all $j$ and all $z\in \partial G_\epsilon\subset (\bc^*)^n$. For the second inequality, we calculate
\[
\det(H^{2\times 2}_j) = (x_j^2+y_j^2)^2(\log|z_j|-(\log|z_j|)^2),
\]
hence
\al{
\det(H^{2\times 2}_j) >0 &\Longleftrightarrow~ \log|z_j|-(\log|z_j|)^2>0\\
&\Longleftrightarrow~ \log|z_j|\in(0,1)\\
&\Longleftrightarrow~ |z_j|\in(1,e).
}
So at $p=(e^{\sqrt{\epsilon/n}},\dots,e^{\sqrt{\epsilon/n}})\in\partial G_\epsilon$ we have $\det\big(H^{2\times 2}_j(p)\big)>0$ for each $j=1,\dots,n$. We conclude that $H\rho(p)$ is positive definite. \qed

\begin{Rmk} \em{The condition $\epsilon<n$ is necessary and sufficient for the existence of some point in $\partial G_\epsilon$ where the real Hessian $H\rho$ is positive definite. This is because when $\epsilon\geq n$, for any $p'=(z_1,\dots,z_n)\in \partial G_\epsilon$ there exists at least one $j=1,\dots,n$ such that $(\log|z_j|)^2\geq 1$. Thus $\det\big(H^{2\times 2}_j(p')\big)\leq 0$, hence $H\rho(p')$ is not positive definite.}
\end{Rmk}

At the point $p=(e^{\sqrt{\epsilon/n}},\dots,e^{\sqrt{\epsilon/n}})\in\partial G_\epsilon$, since $H\rho(p)$ is positive definite and $\rho$ is smooth, there exists some open convex neighborhood $U_p$ of $p$ in $(\bc^*)^n$ (e.g. a sufficiently small open ball centered at $p$), such that $H\rho(z)$ is positive definite for all $z\in U_p$. Thus $\rho$ is strictly convex in $U_p$, hence for any $p_1,p_2\in G_\epsilon\cap U_p$ and any $t\in(0,1)$ we have
\[
\rho\big(tp_1+(1-t)p_2\big)<t\rho(p_1)+(1-t)\rho(p_2)<\epsilon.
\]
So $tp_1+(1-t)p_2\in G_\epsilon$ and is also contained in $U_p$ since $U_p$ is convex. So $tp_1+(1-t)p_2\in G_\epsilon\cap U_p$. We proved that $G_\epsilon\cap U_p$ is convex.

Now we define $\varOmega:=B(2\sqrt{n}e^{\sqrt{\epsilon}})$, $\varOmega':=(\varOmega\backslash \overline{G_\epsilon})\cup U_p$ and $K:=\varOmega\backslash (G_\epsilon\cup \varOmega')$. In fact we have $K=\partial G_\epsilon\backslash U_p\subset\subset \varOmega$.

\begin{Prop}\label{prop-omega} The open set $\varOmega'\cup G_\epsilon=\varOmega\backslash K$ is connected.
\end{Prop}

To prove the proposition, we will use the language in Range's book \cite{Range} Chap 3.7. We call a compact set $A\subset \bc^n$ a {\sl Stein compactum} if it has a neighborhood basis of Stein domains. 

\begin{Lem}\label{lem:Steincompactum} For any bounded Stein compactum $A\subset \bc^n$ $(n\geqslant 2)$, the complement $A^c:=\bc^n\backslash A$ is connected.
\end{Lem}

\proof Since $A$ is bounded, there exists some $R>0$ such that $A\subset B(R)\subset \bc^n$. Thus the open set $A^c$ has an unbounded component containing $B(R)^c$.

Suppose $A^c$ is not connected and $W$ is another component, then $W\subset B(R)$ is bounded. Also $\partial W\subset\partial A$. Take $p\in W\subset A^c$. We have $p\notin A$. Thus $\bc^n\backslash\{p\}$ is an open neighborhood of $A$. Since $A$ is a Stein compactum, there exists an open Stein neighborhood $V$ of $A$ contained in $\bc^n\backslash\{p\}$. The set $\partial W\subset \partial A\subset A$ is relatively compact in $V$. The holomorphically convex hull of $\partial W$ with respect to $V$
\al{
\widehat{\partial W}_V:=\{z\in V||f(z)|\leqslant \sup\limits_{w\in \partial W}|f(w)|,\forall f\in\co(V)\},
}
on one hand, should be relatively compact in $V$ since $V$ is Stein.

But on the other hand $W\backslash V\subset\subset V\cup W$ and $(V\cup W)\backslash (W\backslash V)=V$ is connected. Thus by Hartogs' extension theorem for holomorphic functions we have $\co(V)=\co(V\cup W)$ and by maximal principle we have
\al{
|f(z)|\leqslant\sup\limits_{w\in\partial W}|f(w)|,\cond{z\in W, f\in \co(V\cup W)\subset\co(W)}.
}
Hence $V\cap W\subset\widehat{\partial W}_V$. If $\widehat{\partial W}_V$ is relatively compact in $V$, so is $W\cap V$. But since $W$ is a connected domain, it is pathly connected. We take a point $q\in W\cap V$ and a path $\gamma$ in $W$ connecting $q$ and $p\in W\backslash V$. We get a subpath in $W\cap V$ approaching $\partial V$. Hence $W\cap V$ is not relatively compact in $V$, a contradiction.
\qed

\proof[Proof of Proposition \ref{prop-omega}] In fact $U_p$ is a domain meeting both $G_\epsilon$ and $\varOmega\backslash\overline{G_\epsilon}$. By Proposition \ref{prop:G} {\bf(iii)} we know $G_\epsilon$ is connected. So it suffices to show that $\varOmega\backslash\overline{G_\epsilon}$ is connected. Note that the bounded compact set $\overline{G_\epsilon}$ has a neighborhood basis of Stein domains $\{G_{\epsilon+\delta},\delta>0\}$. By Lemma \ref{lem:Steincompactum} we know $\overline{G_{\epsilon}}^c$ is connected. Since $\overline{G_{\epsilon}}\subset\subset\varOmega=B(2\sqrt{n})$, we know $\varOmega\backslash\overline{G_\epsilon}$ is connected.\qed

\proof [End of the proof of Theorem \ref{main-counterexample}]
Take the trivial holomorphic line bundle $L_{\text{triv}}$ over $\varOmega'$, take $L_{\text{nt}}$ a nontrivial holomorphic line bundle over $G_\epsilon$. Since $G_\epsilon\cap U_p$ is convex, by Corollary \ref{crit-2} we know the restrictions of $L_{\text{nt}}$ and $L_{\text{triv}}$ to $G_\epsilon\cap U_p$ are trivial, hence isomorphic. By the gluing Lemma \ref{lem:gluing} we get a holomorphic line bundle $L_{\text{cex}}$ over $G_\epsilon\cup \varOmega'=\varOmega\backslash K$, which is nontrivial since $L_{\text{cex}}|_{G_\epsilon}\cong L_{\text{nt}}$ is. Thus $H^1(\varOmega\backslash K,\co^*)\neq 0$. However, $H^1(\varOmega,\co^*)=0$ since $\varOmega$ is convex. Thus the restriction map
\al{H^1(\varOmega,\co^*)\rw H^1(\varOmega\backslash K,\co^*)}
cannot be surjective. In particular, $L_{\text{cex}}$ cannot be extended to $\varOmega$.
\endproof

\begin{center}
\begin{picture}(0,0)%
\includegraphics{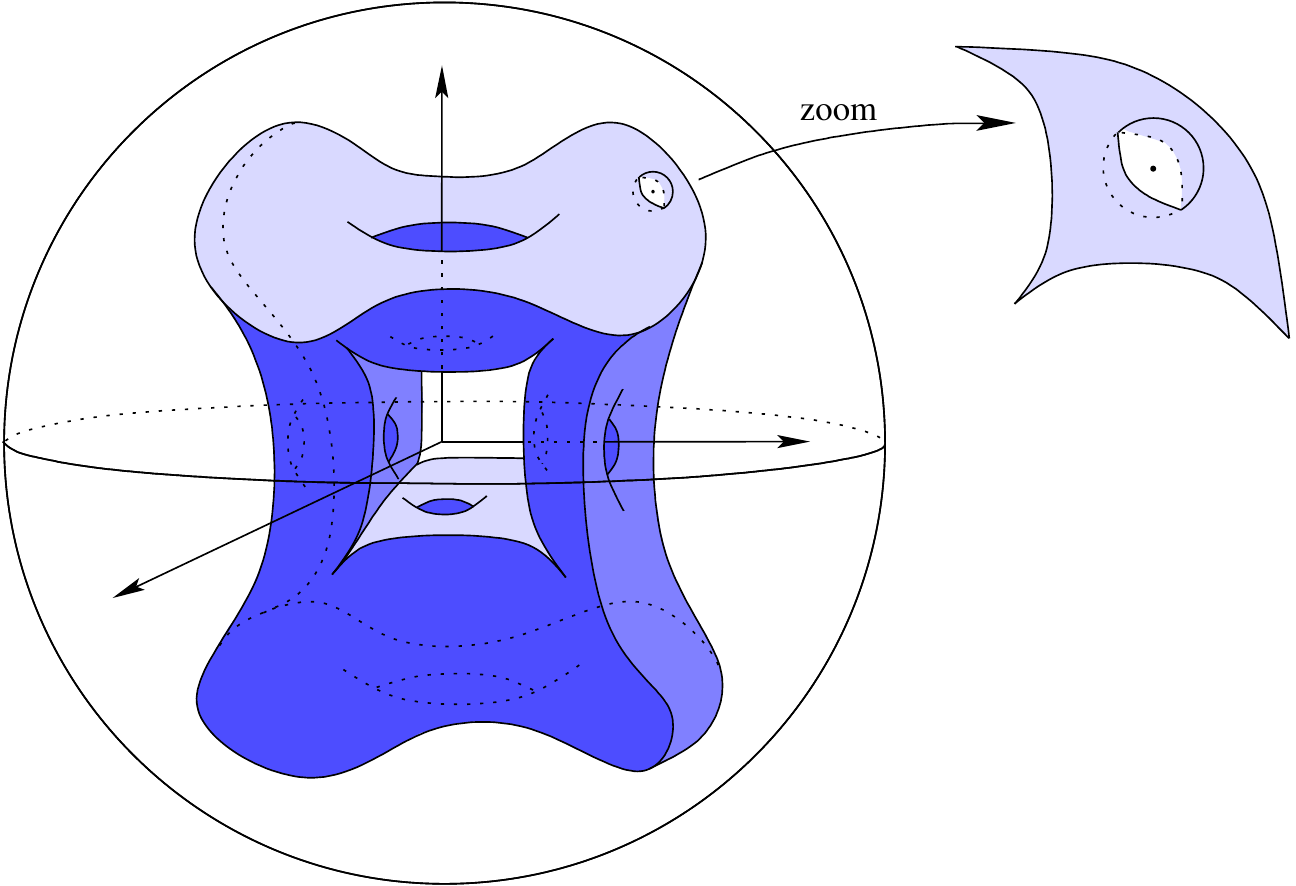}%
\end{picture}%
\setlength{\unitlength}{4144sp}%
\begingroup\makeatletter\ifx\SetFigFont\undefined%
\gdef\SetFigFont#1#2#3#4#5{%
  \reset@font\fontsize{#1}{#2pt}%
  \fontfamily{#3}\fontseries{#4}\fontshape{#5}%
  \selectfont}%
\fi\endgroup%
\begin{picture}(5906,4044)(-11,-3198)
\put(5287,-25){\makebox(0,0)[lb]{\smash{{\SetFigFont{10}{12.0}{\familydefault}{\mddefault}{\updefault}{\color[rgb]{0,0,0}$p$}%
}}}}
\put(4996,344){\makebox(0,0)[lb]{\smash{{\SetFigFont{10}{12.0}{\familydefault}{\mddefault}{\updefault}{\color[rgb]{0,0,0}$U_p$}%
}}}}
\put(1126,-2401){\makebox(0,0)[lb]{\smash{{\SetFigFont{10}{12.0}{\familydefault}{\mddefault}{\updefault}{\color[rgb]{0,0,0}$G_\epsilon$}%
}}}}
\put(1081,-97){\makebox(0,0)[lb]{\smash{{\SetFigFont{10}{12.0}{\familydefault}{\mddefault}{\updefault}{\color[rgb]{0,0,0}$K=\partial G_\epsilon\backslash U_p$}%
}}}}
\put(1351,587){\makebox(0,0)[lb]{\smash{{\SetFigFont{10}{12.0}{\familydefault}{\mddefault}{\updefault}{\color[rgb]{0,0,0}$\varOmega=B(2\sqrt{n}e^{\sqrt{\epsilon}})$}%
}}}}
\end{picture}%

\end{center}

\bibliography{ext}
\bibliographystyle{plain}

\end{document}